\documentclass[12pt]{article}
\usepackage{amssymb}
\usepackage{amsmath}
\usepackage{amsthm}
\usepackage{color}
\usepackage{comment}
\usepackage{geometry}
\usepackage{graphicx}
\usepackage{authblk}

\makeatletter
	\@addtoreset{equation}{section}
\makeatother

\let\originalleft\left
\let\originalright\right
\renewcommand{\left}{\mathopen{}\mathclose\bgroup\originalleft}
\renewcommand{\right}{\aftergroup\egroup\originalright}

\begin{document}

\newcommand{\bO}{{\bf 0}}
\newcommand{\bz}{{\bf z}}
\newcommand{\cA}{\mathcal{A}}
\newcommand{\cN}{\mathcal{N}}
\newcommand{\cO}{\mathcal{O}}
\newcommand{\cP}{\mathcal{P}}
\newcommand{\cT}{\mathcal{T}}
\newcommand{\omDN}{\omega_1}
\newcommand{\omF}{\omega}
\newcommand{\rD}{{\rm D}}
\newcommand{\re}{{\rm e}}
\newcommand{\ri}{{\rm i}}
\newcommand{\myStep}[2]{{\bf Step #1} --- #2\\}

\newtheorem{theorem}{Theorem}[section]

\title{
The VIVID function for numerically continuing periodic orbits arising from grazing bifurcations of hybrid dynamical systems.
}
\author{
I.~Ghosh$^1$\footnote{Author to whom correspondence must be addressed (indra.ghosh@ucd.ie).} and D.J.W.~Simpson$^2$\\

$^1$ School of Mathematics and Statistics,
University College Dublin,
Dublin, D04 V1W8, Ireland\\

$^2$ School of Mathematical and Computational Sciences,
Massey University, 
Palmerston North, 4410,
New Zealand
}

\maketitle



\begin{abstract}

Periodic orbits of systems of ordinary differential equations
can be found and continued numerically by following fixed points of Poincar\'e maps.
However, this often fails near grazing bifurcations where a periodic orbit
collides tangentially with a boundary of phase space.
Failure occurs when the map contains a square-root singularity
and the root-finding algorithm searches beyond the domain of viable values.
We show that by instead following the zeros of a function
that maps Velocity Into Variation In Displacement (VIVID)
this issue is circumvented and there is no such failure.
We illustrate this with a prototypical one-degree-of-freedom impact oscillator model
by applying Newton's method to the VIVID function to follow periodic orbits collapsing into grazing bifurcations.
We also follow curves of saddle-node and period-doubling bifurcations of periodic orbits
that issue from a codimension-two resonant grazing bifurcation.
The VIVID function provides a simple alternative to the more sophisticated collocation method
and enables periodic orbits and their bifurcations to be resolved easily and accurately near grazing bifurcations.

\end{abstract}

\section{Introduction}
\label{sec:intro}

Many mechanical apparatus contain moving parts that bump into one another.
Examples include gear assemblies \cite{HaWi07,ThNa00},
heat exchangers \cite{DeFr99,PaLi92},
atomic force microscopes \cite{DaZh07,MiDa10},
church bells \cite{BrCh18},
and the internal elements of self-propelling capsules for endoscopy \cite{LiPa20,LiWi13}.
These examples involve near misses as well as impacts, so are directly relevant to the material in this paper.

To predict the motions of such apparatus,
one often uses a {\em hybrid model} that combines differential equations and maps \cite{Lu12,VaSc00}. 
The differential equations capture the motion between impacts, while the maps encode the effects of the impacts \cite{AwLa03,BlCz99,Ib09}.
Each map is applied when the system state reaches the corresponding {\em impacting surface} in the phase space of the system.

To understand of the typical long-term dynamics,
and how this differs for different values of the parameters,
we study bifurcations where the dynamics changes in a fundamental way.
In particular, {\em grazing bifurcations} often correspond to the onset of recurring impacts,
and thus are of central importance to many hybrid systems.
A stable non-impacting periodic orbit undergoes a grazing bifurcation
when it first encounters an impacting surface, necessarily intersecting it tangentially.
The bifurcation usually causes the periodic orbit to be destroyed,
being replaced by high-period or chaotic motion.

Since impacts are highly destabilising,
any stable periodic orbits existing near the bifurcation and involving impacts,
usually involve only one impact per period.
Such orbits are termed {\em maximal} \cite{ChOt94}.
They involve some $p \ge 1$ loops close to the path of the non-impacting periodic orbit,
where one loop involves an impact and the remaining $p-1$ loops are near misses.
As a maximal periodic orbit emanates from the grazing bifurcation,
it is typically unstable but may quickly gain stability through either
a saddle-node or period-doubling bifurcation \cite{Iv93,Fo94,Pe96,No01}.

To numerically find and continue periodic orbits, both stable and unstable, two methods are commonly used.
The {\em collocation method} approximates the periodic orbit with an ordered collection of $N$ collocation points
(often with $N$ in the order of $100$) \cite{GoKu05b}.
These points are adjusted iteratively, e.g.~via Newton's method,
until curves through small subsets of consecutive points obey the differential equations to within some tolerance.
This method is effective for smooth differential equations and
has long been embedded into continuation software including {\sc auto} \cite{DoCh07} and {\sc xppaut} \cite{Er02}.
It can be modified to accommodate hybrid systems
by imposing constraints, e.g.~that one collocation point belongs
to an impacting surface and another collocation point is the image of this point under the map.
This is implemented in the software {\sc coco} \cite{DaSc13},
and we believe that previous works where maximal periodic orbits were continued numerically
have exclusively used this approach \cite{DaZh05,JiCh17,MaPi09,ThZh06,YiWe20,ZhDa06}.

The {\em shooting method} instead searches for fixed points of an iterate of a Poincar\'e map.
This is effective for smooth differential equations because the induced Poincar\'e map is smooth (locally),
and fixed points are readily located through root-finding.
In our setting the required Poincar\'e map $P$ is piecewise-smooth,
having one piece $P_{\rm global}$ corresponding to trajectories that miss the impacting surface,
and a second piece $P_{\rm global} \circ P_{{\rm disc},R}$ corresponding to trajectories that experience an impact.
A $p$-loop maximal periodic orbit corresponds to a fixed point of $P_{\rm global}^p \circ P_{{\rm disc},R}$.
This composition is smooth, therefore standard root-finding algorithms can, in principle, be employed to compute the fixed point.
However, near grazing the root-finding algorithm often attempts to evaluate the composition
at points outside its domain of definition, and consequently is unable to locate the fixed point.
If the components of $P$ can be smoothly extended through the switching manifold of $P$,
as is the case for grazing-sliding bifurcations of Filippov systems \cite{DiKo02},
this problem can be alleviated by replacing the components of $P$ with their smooth extensions.
But for typical grazing bifurcations of impacting hybrid systems,
$P_{{\rm disc},R}$ has a square-root singularity and hence does not admit a smooth extension. 
Essentially the algorithm `falls off' the end of the square-root and terminates unsuccessfully.

The purpose of this paper is to overcome this problem by introducing
a new function $V$ that is smooth, has no square-root singularity,
and whose zeros match one-to-one to the fixed points of $P_{\rm global}^p \circ P_{{\rm disc},R}$.
In the context of simple impacting systems,
$V$ maps the velocity of the object at impact
to the variation in displacement of the object following $p$ loops with one impact.
Consequently we call $V$ the VIVID function (abbreviating Velocity Into Variation In Displacement).

Below we illustrate the effectiveness of the VIVID function with a harmonically forced linear oscillator with hard impacts.
We use parameter values based on the physical experiments described in Pavlovskaia {\em et al.}~\cite{PaIn10}
for which a stable two-loop maximal periodic orbit was observed shortly after the grazing bifurcation.
We compute this orbit near grazing and identify saddle-node and period-doubling bifurcations where it gains stability.
The corresponding zero of $V$ is readily continued into, and even through, the grazing bifurcation.
On the other side of the grazing bifurcation the zero does not correspond to a valid periodic orbit
of the hybrid system (here the orbit is {\em virtual} \cite{DiBu08}),
but this means that the solution can be continued as close to the grazing bifurcation
as desired without loss of accuracy or needing to reduce the step size.
This allows us to numerically continue the saddle-node and period-doubling bifurcations, as a second parameter is varied,
until these bifurcations terminate at a codimension-two resonant grazing bifurcation.

The remainder of this paper is organised as follows.
In \S\ref{sec:construct} we define the VIVID function for general hybrid systems with one impacting surface.
We also show how its derivatives can be evaluated for use with Newton's method.
In \S\ref{sec:oscillator} we introduce the linear impact oscillator and identify analytically the grazing bifurcation of the non-impacting periodic solution.
We also derive formulas for the derivatives of the VIVID function in terms of the derivatives of the flow between impacts.
A closed-form expression for the flow is available, thus numerical approximations to these derivatives are not needed,
but numerical root-finding is still required to obtain the times at which orbits experience impacts.
In \S\ref{sec:bifurcations} we showcase one and two parameter bifurcation diagrams of the impact oscillator model produced using the VIVID function.
Finally in \S\ref{sec:conc} we provide concluding remarks.

\section{The definition of the VIVID function}
\label{sec:construct}

In this section we precisely define the VIVID function.
First in \S\ref{sub:generalSystem} we introduce a general class of impacting hybrid systems with one impacting surface.
Then in \S\ref{sub:PGlobalPDisc} we construct a Poincar\'e map $P$
using the discontinuity map framework to accommodate impacts and grazing.
Then in \S\ref{sub:VIVID} we compose the smooth components of this map in an alternate order to build the VIVID function,
and derive formulas for its derivatives.

\subsection{Impacting hybrid systems}
\label{sub:generalSystem}

Let $x(t) \le 0$ denote the position of an object relative to a fixed wall,
and let $y(t) = \dot{x}(t)$ denote its velocity.
We assume that the motion of the object obeys an $n$-dimensional hybrid system of the form
\begin{equation}
\begin{split}
(\dot{x},\dot{y},\dot{\bz}) &= f(x,y,\bz;\mu), \quad \text{for $x < 0$}, \\
(y,\bz) &\mapsto \Phi(y,\bz;\mu), \quad \text{whenever $x = 0$},
\end{split}
\label{eq:hybridSystem}
\end{equation}
where $\bz \in \mathbb{R}^{n-2}$ represents all other state variables and $\mu \in \mathbb{R}$ is a parameter.
The system is written in autonomous form, meaning that any time-dependency (e.g.~forcing) has been incorporated into $\bz$.
While $x(t) < 0$ the system state evolves under the vector field $f$,
and whenever $x(t) = 0$ the system state maps instantaneously under the reset law $\Phi$.
Since $y(t) = \dot{x}(t)$, 
the first component of $f$ is
\begin{equation}
f_1(x,y,\bz;\mu) = y.
\label{eq:f1}
\end{equation}

Let $\Sigma$ denote the coordinate plane $x=0$ in phase space.
This plane is the {\em impacting surface} of \eqref{eq:hybridSystem}.
In view of \eqref{eq:f1}, trajectories following $f$ arrive at $\Sigma$ at points with $y > 0$,
and depart from $\Sigma$ at points with $y < 0$.
For this reason we call the part of $\Sigma$ with $y > 0$ the {\em incoming set},
the part of $\Sigma$ with $y < 0$ the {\em outgoing set},
and the part of $\Sigma$ with $y = 0$ the {\em grazing set}.

We assume that the reset law $\Phi$ has no effect on grazing trajectories.
Thus
\begin{equation}
\Phi(y,\bz;\mu) = (y,\bz) + y \Psi(y,\bz;\mu),
\label{eq:Phi}
\end{equation}
for a smooth function $\Psi$.
We further assume that $\Phi$ causes the object to reverse velocity, possibly with energy loss,
hence maps the incoming set to the outgoing set.

\subsection{Global and discontinuity maps for grazing}
\label{sub:PGlobalPDisc}

Now suppose that when $\mu = 0$ the system
has a periodic orbit passing through $(x,y,\bz) = (0,0,\bz_{\rm graz})$,
for some $\bz_{\rm graz} \in \mathbb{R}^{n-2}$,
but otherwise lying entirely within $x < 0$.
The condition
\begin{equation}
f_2(0,0,\bz_{\rm graz};0) < 0
\label{eq:f2}
\end{equation}
ensures that this orbit
has a quadratic tangency with $\Sigma$ at $(x,y,\bz) = (0,0,\bz_{\rm graz})$,
as is generically the case.

Let $\Pi$ denote the coordinate plane $y=0$ in $\mathbb{R}^n$. 
We use $\Pi$ to define a Poincar\'e map $P : \Pi \to \Pi$ to capture the dynamics near the grazing trajectory,
with the standard discontinuity map modification \cite{DiBu08,FrNo97,FrNo00,No91} to account for impact events.
This assumes $f$ can be smoothly extended beyond $x = 0$ into $x > 0$, as is the case for typical models.

\begin{figure}[b!]
\centering
\includegraphics[width=12cm]{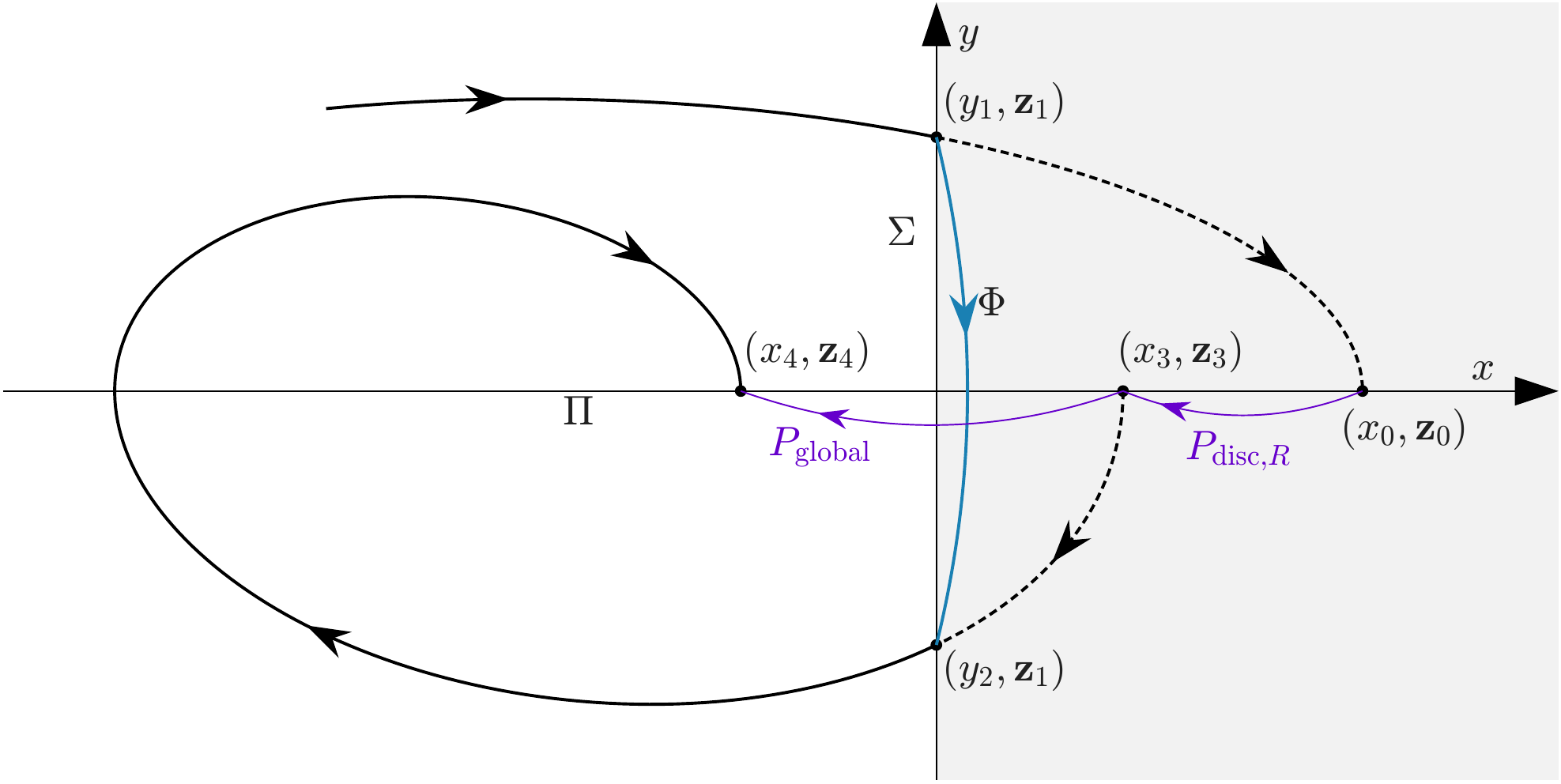}
\caption{
A typical trajectory of the hybrid system \eqref{eq:hybridSystem}
including its extensions (dashed) following $f$ into $x > 0$.
The trajectory intersects the impacting surface $\Sigma$ at $(y_1,\bz_1)$,
then its extension intersects the Poincar\'e section $\Pi$ at $(x_0,\bz_0) = P_{\rm virt}(y_1,\bz_1)$.
The impact point $(y_1,\bz_1)$ maps under the reset law $\Phi$ to $(y_2,\bz_2)$,
then the (backwards) extension of the trajectory from $(y_2,\bz_2)$ intersects $\Pi$ at $(x_3,\bz_3) = P_{\rm virt}(y_2,\bz_2)$.
The map from $(x_0,\bz_0)$ to $(x_3,\bz_3)$
is the right piece $P_{{\rm disc},R}$ of the discontinuity map $P_{\rm disc}$.
The map from $(x_3,\bz_3)$ to the next intersection $(x_4,\bz_4)$ of the trajectory (including its extensions)
with $\Pi$ after one loop following close to the grazing trajectory is the global map $P_{\rm global}$.
\label{fig:PoincareConstruction}}
\end{figure}

We write the Poincar\'e map as $P = P_{\rm global} \circ P_{\rm disc}$,
where $P_{\rm global} : \Pi \to \Pi$ is the {\rm global map} defined by following the flow induced by $f$ forwards in time,
and $P_{\rm disc}$ is the {\em discontinuity map} that accounts for impact events, see Fig.~\ref{fig:PoincareConstruction}.
At points with $x \le 0$, the discontinuity map is the identity map, and so $P_{\rm disc}$ has the piecewise form
\begin{equation}
P_{\rm disc}(x,\bz;\mu) = \begin{cases}
(x,\bz), & x \le 0, \\
P_{{\rm disc},R}(x,\bz;\mu), & x \ge 0.
\end{cases}
\label{eq:Pdisc}
\end{equation}
This is because any point $(x_0,\bz_0)$ in the domain of $P$
is the result of following a trajectory of $f$ until reaching $\Pi$.
If $x_0 < 0$, then no impact has occurred,
while if $x_0 > 0$, then the trajectory of $f$ has passed through $\Sigma$
and the arrival point $(x_0,\bz_0)$ is {\em virtual}.
By evolving forwards under $f$ from $P_{{\rm disc},R}(x_0,\bz_0;\mu) = (x_3,\bz_3)$, see Fig.~\ref{fig:PoincareConstruction},
the subsequent iteration of $P_{\rm global}$ corresponds to the next part of the solution to the full model.

The map $P_{{\rm disc},R}$ can be written as the composition
\begin{equation}
P_{{\rm disc},R} = P_{\rm virt} \circ \Phi \circ P_{\rm virt}^{-1},
\label{eq:PdiscR}
\end{equation}
where $P_{\rm virt} : \Sigma \to \Pi$ is the map defined by following the flow induced by $f$ from $\Sigma$ to $\Pi$.
More specifically, $P_{\rm virt}(y,\bz;\mu)$ is the next point at which the trajectory from $(y,\bz)$ intersects $\Pi$
following $f$ forwards in time if $y > 0$, and backwards in time if $y < 0$, while if $y = 0$ then $P_{\rm virt}$ is the identity map.
The map $P_{\rm virt}$ is well-defined and smooth in a neighbourhood of $(y,\bz;\mu) = (0,\bz_{\rm graz};0)$
because by \eqref{eq:f2} the flow induced by $f$ intersects $\Pi$ transversally.
More specifically, if $f$ is $C^k$, then the flow is $C^k$, and so $P_{\rm virt}$ is $C^k$ as a simple consequence
of the implicit function theorem \cite{Me07}.

For points $(x,z)$ with $x > 0$,
the inverse $P_{\rm virt}^{-1}$ has two branches, one mapping to the incoming set ($y > 0$)
and one mapping to the outgoing set $(y < 0)$.
In \eqref{eq:PdiscR} we use the branch that maps to the incoming set.
This branch is not smooth at $x = 0$ because, as a consequence \eqref{eq:f1} and \eqref{eq:f2},
its leading order contribution is a $\sqrt{x}$-term.

\subsection{Maximal periodic orbits and the VIVID function}
\label{sub:VIVID}

A $p$-loop maximal periodic orbit of \eqref{eq:hybridSystem}
has $p$ loops near the grazing trajectory, exactly one of which experiences an impact.
Let $(y_{\rm imp},\bz_{\rm imp})$ denote the $y$ and $\bz$ values of the orbit at impact,
and let $(x^*,\bz^*) = P_{\rm virt}(y_{\rm imp},\bz_{\rm imp};\mu)$.
Then $(x^*,\bz^*)$ is a fixed point of $Q = P_{\rm global}^p \circ P_{{\rm disc},R}$.
More fully, $Q$ is given by
\begin{equation}
Q = P_{\rm global}^p \circ P_{\rm virt} \circ \Phi \circ P_{\rm virt}^{-1} \,.
\label{eq:Q}
\end{equation}
Since $Q$ is not differentiable at $x = 0$,
we instead consider the {\em VIVID function}
\begin{equation}
V = P_{\rm global}^p \circ P_{\rm virt} \circ \Phi - P_{\rm virt} \,.
\label{eq:V}
\end{equation}
Notice $V(y_{\rm imp},\bz_{\rm imp};\mu) = Q(x^*,\bz^*) - (x^*,z^*) = (0,\bO)$,
so $(y_{\rm imp},\bz_{\rm imp})$ is a zero of $V$.
Indeed, $V = Q \circ P_{\rm virt} - P_{\rm virt}$,
so $(y,\bz)$ is a zero of $V$ if and only if $P_{\rm virt}(y,\bz;\mu)$ is a fixed point of $Q$.

Now assume $\Phi$ can be smoothly extended beyond $y = 0$ into $y < 0$, as is the case for typical models.
Then $V$ is well-defined and smooth on a neighbourhood of $(y,\bz;\mu) = (0,\bz_{\rm graz};0)$
because it is the composition of functions that are each well-defined and smooth locally.
This property is precisely what is needed for numerical continuation, and is not satisfied by $Q$.

To find zeros of $V$ via Newton's method, we need to evaluate $V$ and its derivatives.
Given values of $y_1$, $\bz_1$, and $\mu$, this is achieved as follows.
First we evaluate $(y_2,\bz_2) = \Phi(y_1,\bz_1;\mu)$.
Then we evaluate $(x_0,\bz_0) = P_{\rm virt}(y_1,\bz_1;\mu)$
and $(x_3,\bz_3) = P_{\rm virt}(y_2,\bz_2;\mu)$ by combining a numerical ODE solver (e.g.~Runge-Kutta)
to compute the orbits of $f$ that start from $(0,y_1,\bz_1)$ and $(0,y_2,\bz_2)$,
with a one-dimensional root-finding algorithm (e.g.~the bisection method)
to identify the times at which these orbits intersect $\Pi$.
Then we evaluate $(x_4,\bz_4) = P_{\rm global}^p(y_3,\bz_3;\mu)$
by again using the ODE solver and root-finding to obtain the desired intersection with $\Pi$.
The value of $V$ is then
\begin{equation}
V(y_1,\bz_1;\mu) = \begin{bmatrix} x_4 \\ \bz_4 \end{bmatrix} - \begin{bmatrix} x_1 \\ \bz_1 \end{bmatrix}.
\label{eq:V2}
\end{equation}
By the chain rule, the Jacobian matrix of $V$ and its derivative with respect to $\mu$ are given by
\begin{align}
\rD V(y_1,\bz_1;\mu) &= \rD P_{\rm global}^p(x_3,\bz_3;\mu) \rD P_{\rm virt}(y_2,\bz_2;\mu) \rD \Phi(y_1,\bz_1;\mu)
- \rD P_{\rm virt}(y_1,\bz_1;\mu), \label{eq:dVdyz} \\
\frac{\partial V}{\partial \mu}(y_1,\bz_1;\mu) &=
\rD P_{\rm global}^p(x_3,\bz_3;\mu) \rD P_{\rm virt}(y_2,\bz_2;\mu) \frac{\partial \Phi}{\partial \mu}(y_1,\bz_1;\mu) \nonumber \\
&\quad+ \rD P_{\rm global}^p(x_3,\bz_3;\mu) \frac{\partial P_{\rm virt}}{\partial \mu}(y_2,\bz_2;\mu)
+ \frac{\partial P_{\rm global}^p}{\partial \mu}(x_3,\bz_3;\mu)
- \frac{\partial P_{\rm virt}}{\partial \mu}(y_1,\bz_1;\mu). \label{eq:dVdmu}
\end{align}

In general these derivatives can be evaluated numerically using finite differences.
For the linear impact oscillator studied below,
an explicit expression for the flow of $f$ is available
and so a numerical ODE solver is not needed.
Also we are able to obtain explicit formulas for the components of \eqref{eq:dVdyz} and \eqref{eq:dVdmu},
so finite difference approximations are not required either.
In this case the only numerical approximations used in evaluating $V$ and its first derivatives
are the evolution times obtained by one-dimensional root-finding.

The stability multipliers of the maximal periodic orbit are the eigenvalues of
\begin{equation}
\rD Q(y_1,\bz_1;\mu) = \rD P_{\rm global}^p(x_3,\bz_3;\mu) \rD P_{\rm virt}(y_2,\bz_2;\mu) \rD \Phi(y_1,\bz_1;\mu)
\rD P_{\rm virt}(y_1,\bz_1;\mu)^{-1}.
\label{eq:dQdyz}
\end{equation}
Once a zero of $V$ has been computed, we can evaluate \eqref{eq:dQdyz} and its eigenvalues
to ascertain the stability of the periodic orbit.
Due to the square-root singularity, $\rD P_{\rm virt}(y_1,\bz_1;\mu)$ is singular in the limit $y_1 \to 0$,
and consequently the periodic orbit is typically unstable as it emanates from the grazing bifurcation.

In review, the smoothness of a map defined by following a flow from one surface to another
does not depend on the direction of the flow at the departure surface,
it relies only on the transversality of the flow to the arrival surface and the smoothness of the flow itself.
For a system of the form \eqref{eq:hybridSystem} satisfying \eqref{eq:f1} and \eqref{eq:f2},
the flow of $f$ is transverse to $\Pi$, but not to $\Sigma$.
Thus the map $P_{\rm virt} : \Sigma \to \Pi$ is smooth, but its inverse is not.
By designing the VIVID function in a way that does not require the inverse of $P_{\rm virt}$ to be evaluated,
we have constructed a function that is smooth.

\section{A linear oscillator with hard impacts}
\label{sec:oscillator}

\begin{figure}[b!]
\centering
\includegraphics[width=9cm]{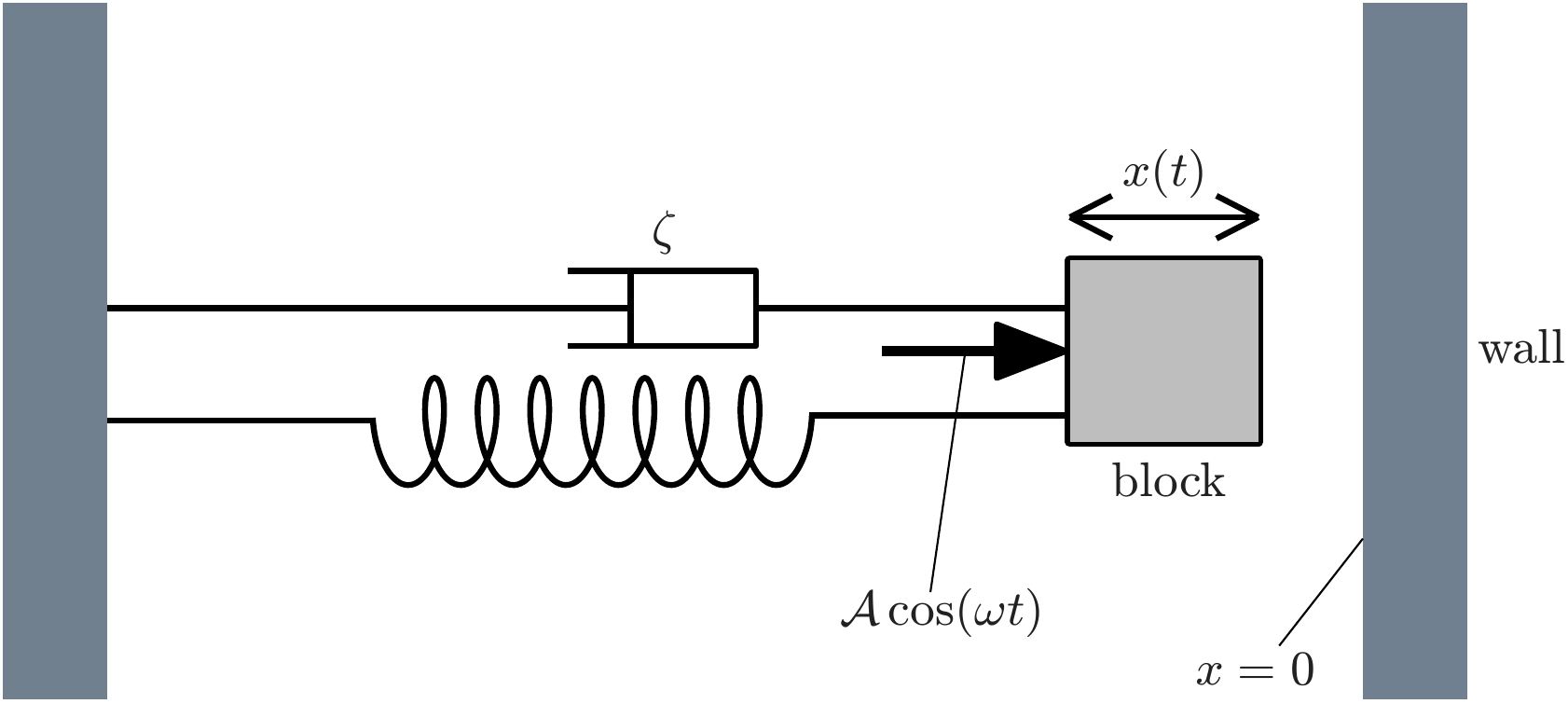}
\caption{
A schematic diagram of the linear impact oscillator modelled by \eqref{eq:impactOsc}.
\label{fig:impactOscillator}}
\end{figure}

As a minimal example of a hybrid system \eqref{eq:hybridSystem},
we consider an underdamped, linear, one-degree-of-freedom oscillator with harmonic forcing
and instantaneous hard impacts with a wall.
The system is illustrated in Fig.~\ref{fig:impactOscillator},
and the non-dimensionalised equations of motion are
\begin{equation}
\begin{split}
\ddot{x} + 2 \zeta \dot{x} + x + 1 &= \cA \cos(\omega t), \quad \text{for $x < 0$}, \\
\dot{x} &\mapsto -\epsilon \dot{x}, \quad \text{whenever $x = 0$},
\end{split}
\label{eq:impactOsc}
\end{equation}
where $\zeta \in (0,1)$ is the damping ratio,
$\cA > 0$ is the forcing amplitude,
$\omega > 0$ is the forcing frequency,
and $\epsilon \in (0,1]$ is the coefficient of restitution.
The equilibrium position of the oscillator is $x = -1$,
and the damped natural frequency of the oscillator is
\begin{equation}
\omega_1 = \sqrt{1 - \zeta^2}.
\label{eq:omega1}
\end{equation}
Studies of \eqref{eq:impactOsc} near grazing
date back to the early 1990's \cite{Fo94,Pe96,No91,FoBi94,LaBu94}
and revealed the dynamics and bifurcation structure to be extraordinarily complex.

In this section we treat the forcing amplitude $\cA$ as the primary bifurcation parameter.
First in \S\ref{sub:grazing} we state an explicit formula $\phi$ for the solution to \eqref{eq:impactOsc} ignoring impacts,
and use this to identify the value $\cA_{\rm graz}$ at which the non-impacting periodic orbit undergoes grazing.
Then in \S\ref{sub:framework} we rewrite \eqref{eq:impactOsc} in the general form \eqref{eq:hybridSystem}.
Finally in \S\ref{sub:formulas} we express the derivatives of the VIVID function in terms of the derivative of $\phi$.

\subsection{Grazing of the non-impacting periodic orbit}
\label{sub:grazing}

Let $\phi(t;x_0,y_0,t_0;\cA)$ denote the solution to the differential equation in \eqref{eq:impactOsc}
with initial condition $(x,\dot{x}) = (x_0,y_0)$ at time $t_0$.
This solution is given by
\begin{equation}
\phi(t;x_0,y_0,t_0;\cA) = \phi_h(t;x_0,y_0,t_0;\cA) + \phi_p(t;\cA),
\label{eq:phi}
\end{equation}
where the particular solution is
\begin{equation}
\phi_p(t;\cA) = -1 + \frac{\cA}{\left( 1-\omega^2 \right)^2 + 4 \zeta^2 \omega^2}
\left( \left( 1-\omega^2 \right) \cos(\omega t) + 2 \zeta \omega \sin(\omega t) \right),
\label{eq:phip}
\end{equation}
and the homogeneous solution is
\begin{align}
\phi_h(t;x_0,y_0,t_0;\cA) &= \re^{-\zeta(t-t_0)} \bigg( \Big( \cos(\omega_1(t-t_0))
+ \frac{\zeta}{\omega_1} \,\sin(\omega_1(t-t_0)) \Big) \big( x_0 - \phi_p(t_0;\cA) \big) \nonumber \\
&\quad+\,\frac{1}{\omega_1} \,\sin(\omega_1(t-t_0)) \Big( y_0 - \frac{d \phi_p}{d t}(t_0;\cA) \Big) \bigg).
\label{eq:phih}
\end{align}
Since the damping ratio $\zeta$ is assumed to be positive,
as $t \to \infty$ the homogeneous solution $\phi_h(t)$ decays to zero,
and thus $\phi(t)$ converges to $\phi_p(t)$.
The particular solution $\phi_p(t)$ is a valid solution
to the hybrid system \eqref{eq:impactOsc} if $\phi_p(t) < 0$ for all $t$.
It is a simple exercise to show from \eqref{eq:phip}
that this occurs when $\cA < \cA_{\rm graz}$, where
\begin{equation}
\cA_{\rm graz} = \sqrt{\left( 1-\omega^2 \right)^2 + 4 \zeta^2 \omega^2}.
\label{eq:Agraz}
\end{equation}
So for $\cA < \cA_{\rm graz}$, $\phi_p(t)$
is a stable non-impacting periodic orbit of \eqref{eq:impactOsc}.
At $\cA = \cA_{\rm graz}$, this orbit grazes the impacting surface $x=0$.

Let
\begin{equation}
z = \omega t ~{\rm mod}~ 2 \pi
\label{eq:z}
\end{equation}
denote the phase of the forcing.
By solving $\phi_p(t) = 0$ with $\cA = \cA_{\rm graz}$,
we find that the phase at grazing is
\begin{equation}
z_{\rm graz} = \begin{cases}
\tan^{-1} \left( \frac{2 \zeta \omega}{1 - \omega^2} \right), & \omega \ne 1, \\
\frac{\pi}{2}, & \omega = 1,
\end{cases}
\label{eq:zGraz}
\end{equation}
with $z_{\rm graz} \in (0,\pi)$.

\subsection{Framework for a local analysis of the grazing bifurcation}
\label{sub:framework}

For the system \eqref{eq:impactOsc}, the vector field $f$
and the reset law $\Phi$ in \eqref{eq:hybridSystem} are given by
\begin{align}
f(x,y,z;\mu) &= \begin{bmatrix} y \\ -x - 2 \zeta y - 1 + (\mu + \cA_{\rm graz}) \cos(z) \\ \omega \end{bmatrix}, &
\Phi(y,z;\mu) &= \begin{bmatrix} -\epsilon y \\ z \end{bmatrix},
\label{eq:fPhi}
\end{align}
where $y = \dot{x}$, $\bz = z$ is a scalar (the phase), and
\begin{equation}
\mu = \cA - \cA_{\rm graz} \,.
\label{eq:mu}
\end{equation}
The genericity condition \eqref{eq:f2} holds because
\begin{equation}
f(0,0,z_{\rm graz};0) = -1 + \cA_{\rm graz} \cos(z_{\rm graz}) = -\omega^2
\label{eq:f2impactOsc}
\end{equation}
is negative.
A local analysis to determine the nature of the invariant sets created at the grazing bifurcation
is not the subject of this paper, see instead \cite{DiBu08,FrNo00},
and in particular Nordmark \cite{No01} for maximal periodic orbits.
For completeness, we provide the formula
\begin{equation}
P_{{\rm disc},R}(x,z;\mu) = \begin{bmatrix}
\epsilon^2 x + \cO \left( \left( \sqrt{x}, z-z_{\rm graz}, \mu \right)^3 \right) \\
-\sqrt{2} (1+\epsilon) \sqrt{x} + z - z_{\rm graz} + \cO \left( \left( \sqrt{x}, z-z_{\rm graz}, \mu \right)^2 \right)
\end{bmatrix},
\label{eq:PdiscRimpactOsc}
\end{equation}
which contains the leading-order terms of the right piece of the discontinuity map
for the system \eqref{eq:impactOsc}, and for a derivation see \cite{DiBu08,FrNo00}.
Notice \eqref{eq:PdiscRimpactOsc} has a $\sqrt{x}$-term in its second component.

\subsection{Formulas for derivatives of the VIVID function}
\label{sub:formulas}

The formulas \eqref{eq:dVdyz} and \eqref{eq:dVdmu}
give the derivatives of the VIVID function in terms of the derivatives of $\Phi$, $P_{\rm virt}$, and $P_{\rm global}^p$.
For the system \eqref{eq:impactOsc},
the derivatives of $\Phi$, $P_{\rm virt}$, and $P_{\rm global}^p$ are given in terms of the derivatives of $\phi$ by
\begin{align}
\rD \Phi(y_1,z_1;\mu) &= \begin{bmatrix}
\frac{\partial y_2}{\partial y_1} & \frac{\partial y_2}{\partial z_1} \\[1.5mm]
\frac{\partial z_2}{\partial y_1} & \frac{\partial z_2}{\partial z_1}
\end{bmatrix}
= \begin{bmatrix}
-\epsilon & 0 \\
0 & 1
\end{bmatrix},
\label{eq:DPhi} \\
\frac{\partial \Phi}{\partial \mu}(y_1,z_1;\mu) &= \begin{bmatrix}
\frac{\partial y_2}{\partial \mu} \\[1.5mm]
\frac{\partial z_2}{\partial \mu}
\end{bmatrix}
= \begin{bmatrix} 0 \\ 0 \end{bmatrix},
\label{eq:dPhidmu} \\
\rD P_{\rm virt}(y_2,z_2;\mu) &= \begin{bmatrix}
\frac{\partial x_3}{\partial y_2} & \frac{\partial x_3}{\partial z_2} \\[1.5mm]
\frac{\partial z_3}{\partial y_2} & \frac{\partial z_3}{\partial z_2}
\end{bmatrix}
= \frac{1}{\frac{\partial^2 \phi}{\partial t^2}}
\begin{bmatrix}
\frac{\partial \phi}{\partial y_0} \frac{\partial^2 \phi}{\partial t^2}
- \frac{\partial \phi}{\partial t} \frac{\partial^2 \phi}{\partial y_0 \partial t} &
\frac{1}{\omega} \left( \frac{\partial \phi}{\partial t_0} \frac{\partial^2 \phi}{\partial t^2}
- \frac{\partial \phi}{\partial t} \frac{\partial^2 \phi}{\partial t_0 \partial t} \right) \\
-\omega \frac{\partial^2 \phi}{\partial y_0 \partial t} &
-\frac{\partial^2 \phi}{\partial t_0 \partial t}
\end{bmatrix}
\Bigg|_{(t_3;0,y_2,t_2;\cA)} \,,
\label{eq:DPvirt} \\
\frac{\partial P_{\rm virt}}{\partial \mu}(y_2,z_2;\mu) &= \begin{bmatrix}
\frac{\partial x_3}{\partial \mu} \\[1.5mm]
\frac{\partial z_3}{\partial \mu}
\end{bmatrix}
= \frac{1}{\frac{\partial^2 \phi}{\partial t^2}}
\begin{bmatrix} \frac{\partial \phi}{\partial \cA} \frac{\partial^2 \phi}{\partial t^2}
- \frac{\partial \phi}{\partial t} \frac{\partial^2 \phi}{\partial \cA \partial t} \\[1.5mm]
-\omega \frac{\partial^2 \phi}{\partial \cA \partial t}
\end{bmatrix}
\Bigg|_{(t_3;0,y_2,t_2;\cA)} \,,
\label{eq:dPvirtdmu} \\
\rD P_{\rm global}^p(x_3,z_3;\mu) &= \begin{bmatrix}
\frac{\partial x_4}{\partial x_3} & \frac{\partial x_4}{\partial z_3} \\[1.5mm]
\frac{\partial z_4}{\partial x_3} & \frac{\partial z_4}{\partial z_3}
\end{bmatrix}
= \frac{1}{\frac{\partial^2 \phi}{\partial t^2}}
\begin{bmatrix}
\frac{\partial \phi}{\partial x_0} \frac{\partial^2 \phi}{\partial t^2}
- \frac{\partial \phi}{\partial t} \frac{\partial^2 \phi}{\partial x_0 \partial t} &
\frac{1}{\omega} \left( \frac{\partial \phi}{\partial t_0} \frac{\partial^2 \phi}{\partial t^2}
- \frac{\partial \phi}{\partial t} \frac{\partial^2 \phi}{\partial t_0 \partial t} \right) \\
-\omega \frac{\partial^2 \phi}{\partial x_0 \partial t} &
-\frac{\partial^2 \phi}{\partial t_0 \partial t}
\end{bmatrix}
\Bigg|_{(t_4;x_3,0,t_3;\cA)} \,,
\label{eq:DPglobalp} \\
\frac{\partial P_{\rm global}^p}{\partial \mu}(x_3,z_3;\mu) &= \begin{bmatrix}
\frac{\partial x_4}{\partial \mu} \\[1.5mm]
\frac{\partial z_4}{\partial \mu}
\end{bmatrix}
= \frac{1}{\frac{\partial^2 \phi}{\partial t^2}}
\begin{bmatrix} \frac{\partial \phi}{\partial \cA} \frac{\partial^2 \phi}{\partial t^2}
- \frac{\partial \phi}{\partial t} \frac{\partial^2 \phi}{\partial \cA \partial t} \\[1.5mm]
-\omega \frac{\partial^2 \phi}{\partial \cA \partial t}
\end{bmatrix}
\Bigg|_{(t_4;x_3,0,t_3;\cA)} \,.
\label{eq:dPglobalpdmu}
\end{align}
In each equation the first equality clarifies the definition of the derivative
by using the points of the trajectory indicated in Fig.~\ref{fig:PoincareConstruction},
while the second equality expresses the derivative in terms of $\phi$ and the parameters of \eqref{eq:impactOsc}.
The formulas \eqref{eq:DPhi} and \eqref{eq:dPhidmu} are trivial consequences of the formula for $\Phi$, \eqref{eq:fPhi}.
The formulas \eqref{eq:DPvirt} and \eqref{eq:dPvirtdmu} are derived in Appendix \ref{app:derivatives}
by calculating an infinitesimally perturbed trajectory,
and \eqref{eq:DPglobalp} and \eqref{eq:dPglobalpdmu} can be derived in the same manner.

In our numerical study of \eqref{eq:impactOsc}
we hard coded formulas for the derivatives of $\phi$
which are available in closed-form in view of \eqref{eq:phi}--\eqref{eq:phih}.
In this way the derivatives of the VIVID function were evaluated without numerical approximation.

\section{Bifurcations and continuation}
\label{sec:bifurcations}

In this section we use the VIVID function
to compute one and two-parameter bifurcation diagrams
of the impact oscillator \eqref{eq:impactOsc}.

\subsection{One-parameter continuation}
\label{sub:one}

We first fix $\zeta = 0.02$, $\epsilon = 0.9$, and $\omega = 0.81$, and vary the forcing amplitude $\cA$.
These values are based on Pavlovskaia {\em et al.}~\cite{PaIn10}, see also \cite{BaIn09,SiAv20},
which compared physical experiments to the numerical simulations of a one-degree-of-freedom impact oscillator model.
Both the experiments and the numerics showed that the stable non-impacting periodic motion
transitioned, through a grazing bifurcation, to chaotic motion over a narrow parameter interval,
and then to a stable two-loop maximal periodic orbit.

\begin{figure}[b!]
\begin{center}
\setlength{\unitlength}{1cm}
\begin{picture}(15.4,6)
\put(.8,0){\includegraphics[height=6cm]{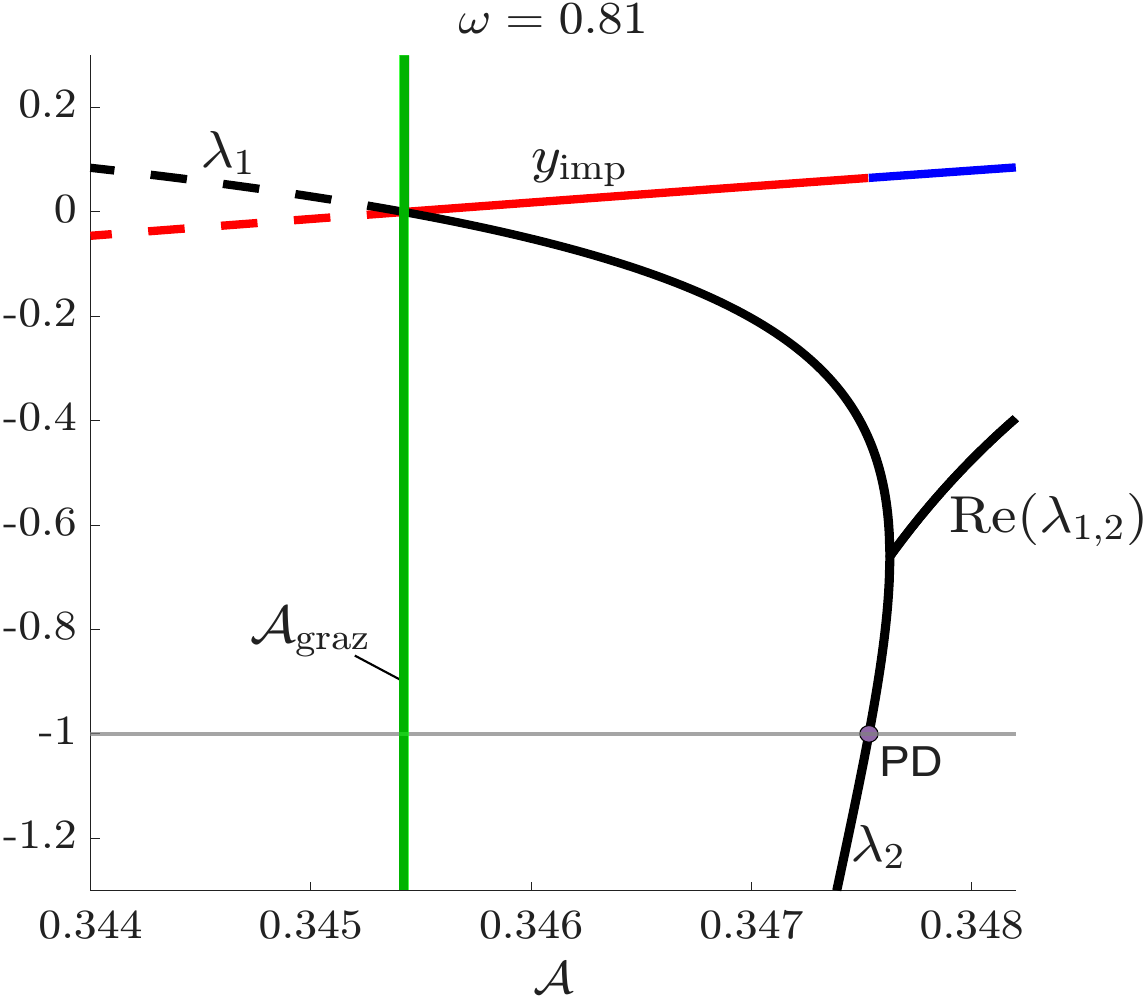}}
\put(8.8,0){\includegraphics[height=6cm]{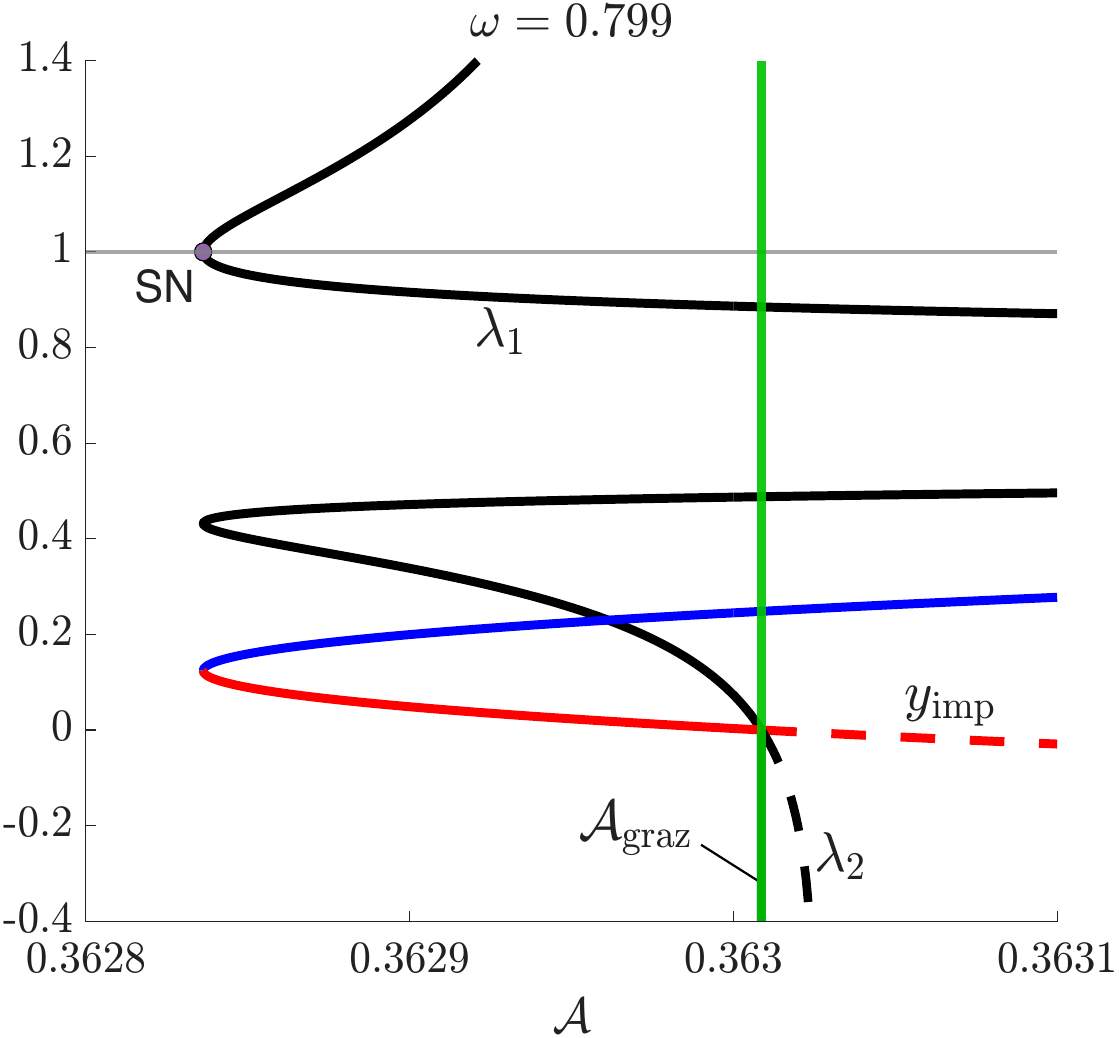}}
\put(0,5.6){\bf \small (a)}
\put(8,5.6){\bf \small (b)}
\end{picture}
\caption{
One-parameter bifurcation diagrams showing the impact velocity $y_{\rm imp}$ and stability multipliers $\lambda_{1,2}$
of a two-loop maximal periodic orbit of \eqref{eq:impactOsc} with $\zeta = 0.02$ and $\epsilon = 0.9$
(PD: period-doubling bifurcation; SN: saddle-node bifurcation).
The curves are blue where the orbit is stable,
red where the orbit is unstable,
and dashed where the orbit is virtual.
The green line indicates the grazing bifurcation of the non-impacting periodic orbit \eqref{eq:phip}.
\label{fig:OneParamBif}}
\end{center}
\end{figure}

Fig.~\ref{fig:OneParamBif}a is a bifurcation diagram of \eqref{eq:impactOsc}
showing how the two-loop maximal periodic orbit varies with $\cA$.
The periodic orbit is unstable (red) as it grows out of the grazing bifurcation $\cA = \cA_{\rm graz}$
and gains stability in a period-doubling bifurcation.
This bifurcation is subcritical, creating an unstable four-loop orbit (not shown).

We now explain how Fig.~\ref{fig:OneParamBif}a was computed.
We began to the right of the period-doubling bifurcation where the two-loop maximal periodic orbit is stable
and approximated it from the long-time behaviour of a numerically computed forward orbit of \eqref{eq:impactOsc}.
To ensure low-velocity impacts were detected
the orbit was computed by iterating the Poincar\'e map $P = P_{\rm global} \circ P_{\rm disc}$
using the bisection method to accurately identify the times at which the orbit intersected $\Sigma$ and $\Pi$.
In particular we obtained values $y_{\rm imp}$ and $z_{\rm imp}$ for the velocity and phase
of the periodic orbit at impact.

We then continued the periodic orbit under variation of the parameter $\cA$.
This can be done by taking steps in $\cA$ and solving for $y_{\rm imp}$ and $z_{\rm imp}$,
but this would not enable us to continue through saddle-node bifurcations (see already Fig.~\ref{fig:OneParamBif}b).
For this reason we took steps in $y_{\rm imp}$ (of size $\Delta y_{\rm imp} = 0.0001$) and solved for $z_{\rm imp}$ and $\cA$.
A more versatile approach is {\em pseudo-arclength continuation} \cite{AlGe03,Do07b,Go00b}
for which no {\em a priori} restrictions are placed upon the variables or parameters at each step,
but steps in $y_{\rm imp}$ are simpler and proved sufficient for our purposes.

At each step $y_{\rm imp}$ is fixed
and initial estimates, $z_0$ and $\cA_0$, for the values of $z_{\rm imp}$ and $\cA$,
were obtained from the outputs of the previous few steps.
The desired values $z_{\rm imp}$ and $\cA$ correspond to a zero of the VIVID function $V$,
using $p=2$ in \eqref{eq:V} for two loops.
These values were obtained through Newton iterations
\begin{equation}
\begin{bmatrix} z_{j+1} \\ \cA_{j+1} \end{bmatrix}
= \begin{bmatrix} z_j \\ \cA_j \end{bmatrix}
- \begin{bmatrix} \frac{\partial v_1}{\partial z} & \frac{\partial v_1}{\partial \mu} \\[1.5mm]
\frac{\partial v_2}{\partial z} & \frac{\partial v_2}{\partial \mu} \end{bmatrix}^{-1}
\begin{bmatrix} v_1 \\ v_2 \end{bmatrix}
\label{eq:Newton}
\end{equation}
where $(v_1,v_2) = V(y_{\rm imp},z_j;\cA_j - \cA_{\rm graz})$.
At each iteration the derivatives in \eqref{eq:Newton} were obtained
by evaluating \eqref{eq:dVdyz} and \eqref{eq:dVdmu}
using the formulas \eqref{eq:DPhi}--\eqref{eq:dPglobalpdmu}.
The iterations were stopped when a tolerance on the norm of $V$ was met.

In this way the two-loop maximal periodic orbit was continued up
to and even through the grazing bifurcation at $\cA = \cA_{\rm graz}$
where the two loops merge and the orbit becomes coincident to the single-loop grazing trajectory.
Beyond the grazing bifurcation the periodic orbit
is not a viable solution to \eqref{eq:impactOsc} as it involves a negative impact velocity
and shown dashed in Fig.~\ref{fig:OneParamBif}a.

At each step we also obtained values for the stability multipliers $\lambda_{1,2}$
of the periodic orbit by evaluating the Jacobian matrix \eqref{eq:dQdyz} and computing its eigenvalues.
As shown in Fig.~\ref{fig:OneParamBif}a, as $\cA$ is decreased
the stability multipliers are initially complex-valued,
then become real, then $\lambda_2$ passes through the value $-1$
where the periodic orbit loses stability in a period-doubling bifurcation.
Note that $\lambda_2 \to -\infty$ as $\cA \to \cA_{\rm graz}$
due to the $\sqrt{x}$-term in \eqref{eq:PdiscRimpactOsc}.

Between the grazing bifurcation and period-doubling bifurcation
neither the non-impacting periodic orbit nor the two-loop maximal periodic orbit is stable.
This gap is the basis for the narrow band of chaotic motion reported in \cite{PaIn10}.

Fig.~\ref{fig:OneParamBif}b shows the result of the same numerical procedure using instead $\omega = 0.799$.
In this case the branch of two-loop maximal periodic orbits extends across the grazing bifurcation.
The branch loses stability and reverses direction at a saddle-node bifurcation (where $\lambda_1 = 1$),
then collapses into the grazing bifurcation.
Between the grazing bifurcation and saddle-node bifurcation
the system is bistable, having both a stable non-impacting periodic orbit
and a stable two-loop maximal periodic orbit.

In both examples the grazing bifurcation generates an unstable two-loop maximal periodic orbit.
In Fig.~\ref{fig:OneParamBif}a the orbit emanates to the right of the bifurcation,
whereas in Fig.~\ref{fig:OneParamBif}a the orbit emanates to the left of the bifurcation.
As shown by Nordmark \cite{No01} for arbitrary $p$-loop maximal periodic orbits,
this direction changes at certain rational ratios
between the forcing frequency $\omega$ and the damped natural frequency $\omega_1 = \sqrt{1 - \zeta^2}$.
This is a type of resonance, which here occurs when $\frac{\omega_1}{\omega} = \frac{5}{4}$,
so $\omega = \frac{4}{5} \sqrt{1 - \zeta^2} \approx 0.7998$ using $\zeta = 0.02$.

\subsection{Two-parameter continuation}
\label{sub:two}

\begin{figure}[b!]
\begin{center}
\setlength{\unitlength}{1cm}
\begin{picture}(16,6)
\put(.3,0){\includegraphics[height=6cm]{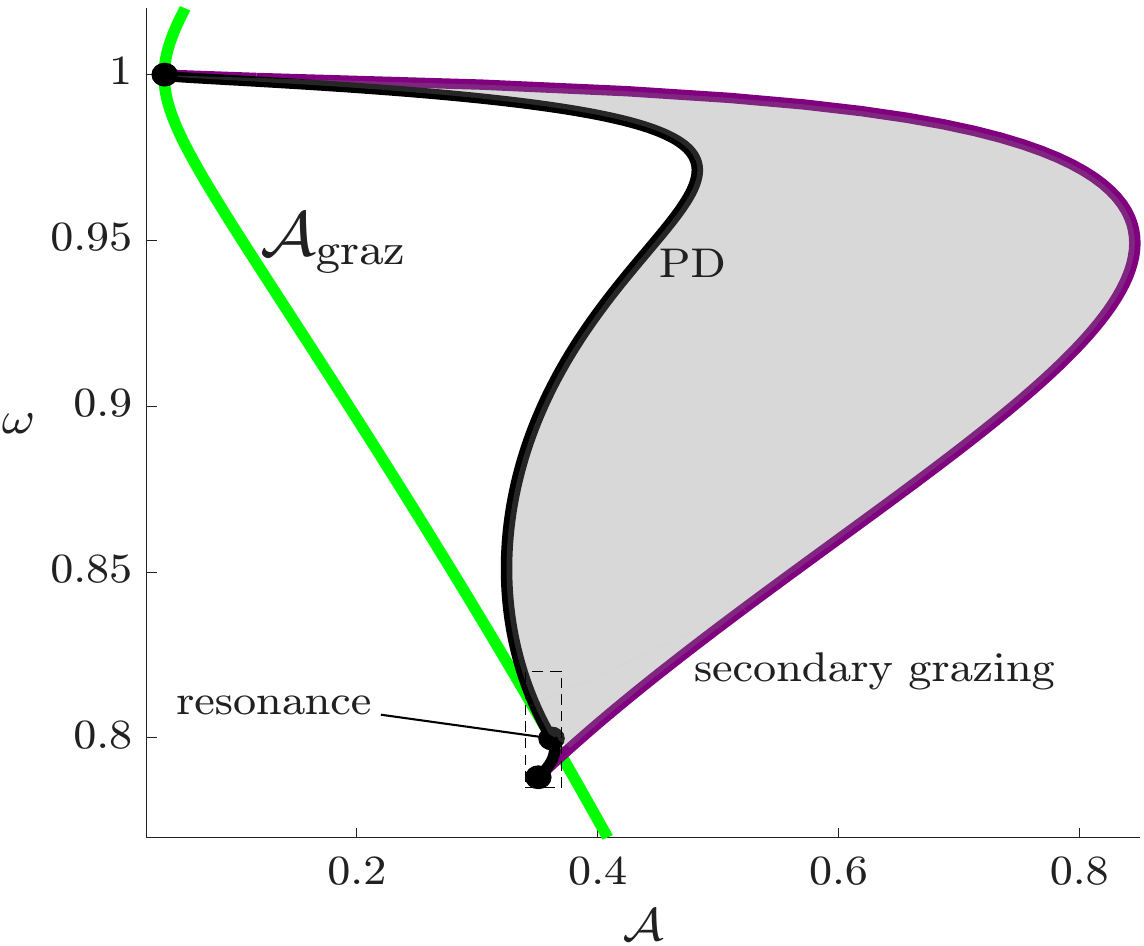}}
\put(9.1,0){\includegraphics[height=6cm]{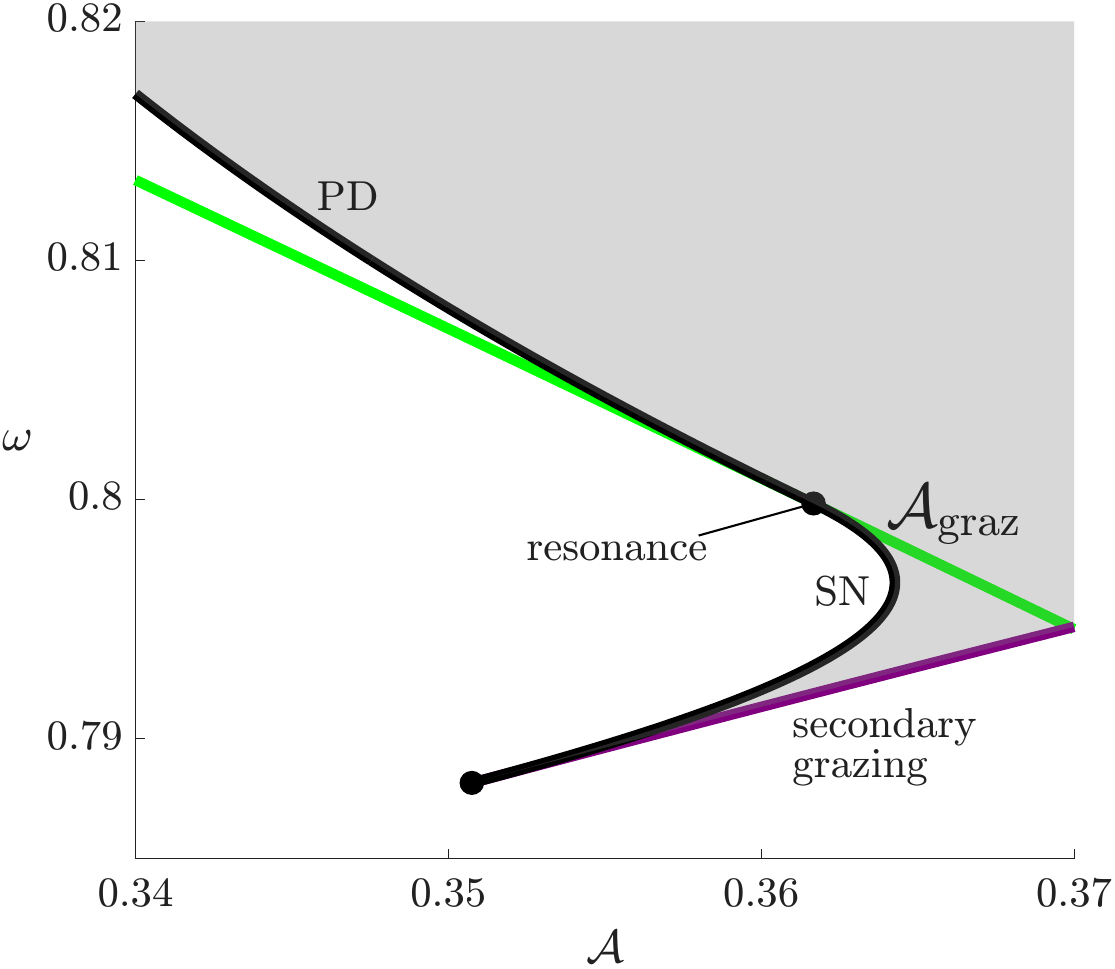}}
\put(0,5.6){\bf \small (a)}
\put(8.5,5.6){\bf \small (b)}
\end{picture}
\caption{
Panel (a) is a two-parameter bifurcation diagram of \eqref{eq:impactOsc}
with $\zeta = 0.02$ and $\epsilon = 0.9$; panel (b) is a magnification of the area within the boundary portrayed by the broken lines in panel (a) 
(PD: period-doubling bifurcation; SN: saddle-node bifurcation).
\label{fig:TwoParamBif}
}
\end{center}
\end{figure}

The codimension-two point of resonance is indicated in Fig.~\ref{fig:TwoParamBif}.
This is a two-parameter bifurcation diagram showing the grazing curve $\cA = \cA_{\rm graz}(\omega)$,
and curves of period-doubling and saddle-node bifurcations of the two-loop maximal periodic orbit.
These curves were continued numerically by performing steps in $\omega$
and using the bisection method to solve $\lambda_2 = -1$ or $\lambda_1 = 1$.
At each step in the continuation, $\omega$ was fixed while $y_{\rm imp}$ was adjusted according to the bisection method.
At each step in the bisection method, $z_{\rm imp}$ and $\cA$ were obtained via the Newton iterations \eqref{eq:Newton}.

The period-doubling and saddle-node curves both terminate at the resonance point.
This phenomenon has been described for similar impacting models by Nordmark \cite{No01}
and Thota {\em et al.}~\cite{ThZh06},
and originally by Ivanov \cite{Iv93} in the case of one-loop maximal periodic orbits.
The other termination point of the period-doubling curve occurs at $\omega \approx 1$
corresponding to the resonance $\frac{\omega_1}{\omega} = 1$.
The other termination point of the saddle-node curve occurs at $\omega \approx 0.7881$
where the two-loop maximal periodic orbit grazes the impacting surface.
From this point there issues a curve of `secondary' grazing bifurcations
of the two-loop maximal periodic orbit, see Fig.~\ref{fig:TwoParamBif}.
This curve, together with the period-doubling and saddle-node curves,
bound the region (shaded) where the system has a stable two-loop maximal periodic orbit.

\section{Discussion}
\label{sec:conc}

In this paper we have introduced the VIVID function
for supporting numerical bifurcation analyses of impacting hybrid systems,
and illustrated its applicability with a one-degree-of-freedom impact oscillator model.
Algebraically the VIVID function \eqref{eq:V} is a minor modification of the usual Poincar\'e map \eqref{eq:Q},
but conceptually it is a major alteration because we are no longer searching for the fixed points of a map,
we are searching for the zeros of a function.

Our numerical results have highlighted the significance of resonant grazing bifurcations.
These are the exceptional grazing bifurcations
at which one of the bifurcating periodic orbits is stable, despite involving impacts.
A general unfolding of these bifurcations remains for future work,
and we believe this can be achieved by employing the VIVID function in a theoretical role.
Since the VIVID function is smooth in a neighbourhood of the grazing bifurcation,
it can be used in conjunction with the implicit function theorem
to formally identify the smooth curves of period-doubling and saddle-node bifurcations.

\section{Data availability statement}
The data that support the findings of this study are openly available at the following URL/DOI: https://github.com/indrag49/VIVID.

\section*{Acknowledgements}

This work was supported by Marsden Fund contract MAU2209 managed by Royal Society Te Ap\={a}rangi.
The authors thank John Bailie and Harry Dankowicz for helpful conversations.

\appendix

\section{Derivation of formulas for the derivatives of the virtual map	}
\label{app:derivatives}

Given $y_2$, $z_2$, and $\cA$, write
\begin{equation}
\begin{bmatrix} x_3 \\ z_3 \end{bmatrix} = P_{\rm virt}(y_2,z_2;\mu),
\nonumber
\end{equation}
where $\mu = \cA - \cA_{\rm graz}$.
The time $t_2$ at which the orbit following $f$ passes through $(x,y,z) = (0,y_2,z_2)$
is $t_2 = \frac{z_2}{\omega} + \frac{2 \pi k}{\omega}$, for some $k \in \mathbb{Z}$,
and without loss of generality we take $k = 0$.
The time $t_3$ at which this orbit passes through $(x,y,z) = (x_3,0,z_3)$
satisfies $\omega t_3 \,{\rm mod}\, 2 \pi = z_3$.

Our task is to evaluate $\rD P_{\rm virt}$ and $\frac{\partial P_{\rm virt}}{\partial \mu}$,
which we write as
\begin{align}
\rD P_{\rm virt}(y_2,z_2;\mu) &= \begin{bmatrix}
a_{11} & a_{12} \\
a_{21} & a_{22}
\end{bmatrix}, &
\frac{\partial P_{\rm virt}}{\partial \mu}(y_2,z_2;\mu) = \begin{bmatrix}
b_1 \\ b_2
\end{bmatrix},
\label{eq:Ab}
\end{align}
and verify \eqref{eq:DPvirt} and \eqref{eq:dPvirtdmu}.

Given small $\delta_1, \delta_2, \delta_3 \in \mathbb{R}$, write
\begin{equation}
\begin{bmatrix} x_3' \\ z_3' \end{bmatrix} = P_{\rm virt}(y_2 + \delta_1,z_2 + \delta_2;\mu + \delta_3).
\nonumber
\end{equation}
Then
\begin{align}
x_3' &= x_3 + a_{11} \delta_1 + a_{12} \delta_2 + b_1 \delta_3 + \cO(2), \nonumber \\
z_3' &= z_3 + a_{21} \delta_1 + a_{22} \delta_2 + b_2 \delta_3 + \cO(2), \nonumber
\end{align}
where $\cO(2)$ denotes all terms that are quadratic or higher order in $\delta_1$, $\delta_2$, and $\delta_3$.

By the definitions of $\phi$ and $P_{\rm virt}$,
\begin{align}
x_3' &= \phi \left( t_3';0,y_2 + \delta_1,t_2 + \tfrac{\delta_2}{\omega};\cA + \delta_3 \right), \label{eq:w1} \\
0 &= \frac{\partial \phi}{\partial t} \left( t_3';0,y_2 + \delta_1,t_2 + \tfrac{\delta_2}{\omega};\cA + \delta_3 \right), \label{eq:w2}
\end{align}
where $t_3' = t_3 + \frac{a_{21} \delta_1}{\omega} + \frac{a_{22} \delta_2}{\omega} + \frac{b_2 \delta_3}{\omega} + \cO(2)$.
By Taylor expanding $\phi$ and $\frac{\partial \phi}{\partial t}$
and retaining only the first-order terms in \eqref{eq:w1} and \eqref{eq:w2}, we obtain
\begin{align}
a_{11} \delta_1 + a_{12} \delta_2 + b_1 \delta_3
&= \frac{1}{\omega} \frac{\partial \phi}{\partial t} \left( a_{21} \delta_1 + a_{22} \delta_2 + b_2 \delta_3 \right)
+ \frac{\partial \phi}{\partial y_0} \delta_1
+ \frac{1}{\omega} \frac{\partial \phi}{\partial t_0} \delta_2
+ \frac{\partial \phi}{\partial \cA} \delta_3, \nonumber \\
0 &= \frac{1}{\omega} \frac{\partial^2 \phi}{\partial t^2} \left( a_{21} \delta_1 + a_{22} \delta_2 + b_2 \delta_3 \right)
+ \frac{\partial^2 \phi}{\partial y_0 \partial t} \delta_1
+ \frac{1}{\omega} \frac{\partial^2 \phi}{\partial t_0 \partial t} \delta_2
+ \frac{\partial^2 \phi}{\partial \cA \partial t} \delta_3, \nonumber
\end{align}
where the derivatives are evaluated at $(t_3;0,y_2,t_2;\cA)$.
By matching the $\delta_1$, $\delta_2$, and $\delta_3$ terms we obtain
\begin{align}
a_{11} &= \frac{1}{\omega} \frac{\partial \phi}{\partial t} a_{21} + \frac{\partial \phi}{\partial y_0}, &
0 &= \frac{1}{\omega} \frac{\partial^2 \phi}{\partial t^2} a_{21} + \frac{\partial^2 \phi}{\partial y_0 \partial t}, \nonumber \\
a_{12} &= \frac{1}{\omega} \frac{\partial \phi}{\partial t} a_{22} + \frac{1}{\omega} \frac{\partial \phi}{\partial t_0}, &
0 &= \frac{1}{\omega} \frac{\partial^2 \phi}{\partial t^2} a_{22} + \frac{1}{\omega} \frac{\partial^2 \phi}{\partial t_0 \partial t}, \nonumber \\
b_1 &= \frac{1}{\omega} \frac{\partial \phi}{\partial t} b_2 + \frac{\partial \phi}{\partial \cA}, &
0 &= \frac{1}{\omega} \frac{\partial^2 \phi}{\partial t^2} b_2 + \frac{\partial^2 \phi}{\partial \cA \partial t}, \nonumber
\end{align}
and by solving these equations simultaneously for $a_{11}$, $a_{12}$, $a_{21}$, $a_{22}$, $b_1$, and $b_2$
and substituting the results into \eqref{eq:Ab} we arrive at \eqref{eq:DPvirt} and \eqref{eq:dPvirtdmu} as required.


\begin{thebibliography}{10}

\bibitem{HaWi07}
C.K. Halse, R.E. Wilson, M.~di~Bernardo, and M.E. Homer.
\newblock Coexisting solutions and bifurcations in mechanical oscillations with
  backlash.
\newblock {\em J. Sound Vib.}, 305:854--885, 2007.

\bibitem{ThNa00}
S.~Theodossiades and S.~Natsiavas.
\newblock Non-linear dynamics of gear-pair systems with periodic stiffness and
  backlash.
\newblock {\em J. Sound Vib.}, 229(2):287--310, 2000.

\bibitem{DeFr99}
J.M. de~Bedout, M.A. Franchek, and A.K. Bajaj.
\newblock Robust control of chaotic vibrations for impacting heat exchanger
  tubes in crossflow.
\newblock {\em J. Sound Vib.}, 227(1):183--204, 1999.

\bibitem{PaLi92}
M.P. Pa\"{\i}doussis and G.X. Li.
\newblock Cross-flow-induced chaotic vibrations of heat-exchanger tubes
  impacting on loose supports.
\newblock {\em J. Sound Vib.}, 152(2):305--326, 1992.

\bibitem{DaZh07}
H.~Dankowicz, X.~Zhao, and S.~Misra.
\newblock Near-grazing dynamics in tapping-mode atomic-force microscopy.
\newblock {\em Int. J. Non-Linear Mech.}, 42(4):697--709, 2007.

\bibitem{MiDa10}
S.~Misra, H.~Dankowicz, and M.R. Paul.
\newblock Degenerate discontinuity-induced bifurcations in tapping-mode.
\newblock {\em Phys. D}, 239:33--43, 2010.

\bibitem{BrCh18}
P.~Brzeski, A.S.E. Chong, M.~Wiercigroch, and P.~Perlikowski.
\newblock Impact adding bifurcation in an autonomous hybrid dynamical model of
  a church bell.
\newblock {\em Mech. Syst. Signal Process.}, 104:716--724, 2018.

\bibitem{LiPa20}
Y.~Liu, J.~P\'{a}ez~Ch\'{a}vez, J.~Zhang, J.~Tian, B.~Guo, and S.~Prasad.
\newblock The vibro-impact capsule system in millimeter scale: numerical
  optimisation and experimental verification.
\newblock {\em Meccanica}, 55:1885--1902, 2020.

\bibitem{LiWi13}
Y.~Liu, M.~Wiercigroch, E.~Pavlovskaia, and H.~Yu.
\newblock Modelling of a vibro-impact capsule system.
\newblock {\em Int. J. Mech. Sci.}, 66:2--11, 2013.

\bibitem{Lu12}
A.C.J. Luo.
\newblock {\em Discontinuous dynamical systems.}
\newblock Springer, New York, 2012.

\bibitem{VaSc00}
A.J. Van~der Schaft and J.M. Schumacher.
\newblock {\em An Introduction to Hybrid Dynamical Systems.}
\newblock Springer-Verlag, New York, 2000.

\bibitem{AwLa03}
J.~Awrejcewicz and C.~Lamarque.
\newblock {\em Bifurcation and Chaos in Nonsmooth Mechanical Systems.}
\newblock World Scientific, Singapore, 2003.

\bibitem{BlCz99}
B.~Blazejczyk-Okolewska, K.~Czolczynski, T.~Kapitaniak, and J.~Wojewoda.
\newblock {\em Chaotic Mechanics in Systems with Impacts and Friction}.
\newblock World Scientific, Singapore, 1999.

\bibitem{Ib09}
R.A. Ibrahim.
\newblock {\em Vibro-Impact Dynamics.}, volume~43 of {\em Lecture Notes in
  Applied and Computational Mechanics.}
\newblock Springer, New York, 2009.

\bibitem{ChOt94}
W.~Chin, E.~Ott, H.E. Nusse, and C.~Grebogi.
\newblock Grazing bifurcations in impact oscillators.
\newblock {\em Phys. Rev. E}, 50(6):4427--4450, 1994.

\bibitem{Iv93}
A.P. Ivanov.
\newblock Stabilization of an impact oscillator near grazing incidence owing to
  resonance.
\newblock {\em J. Sound Vib.}, 162(3):562--565, 1993.

\bibitem{Fo94}
S.~Foale.
\newblock Analytical determination of bifurcations in an impact oscillator.
\newblock {\em Proc. R. Soc. Lond. A}, 347:353--364, 1994.

\bibitem{Pe96}
F.~Peterka.
\newblock Bifurcations and transition phenomena in an impact oscillator.
\newblock {\em Chaos Solitons Fractals}, 7(10):1635--1647, 1996.

\bibitem{No01}
A.B. Nordmark.
\newblock Existence of periodic orbits in grazing bifurcations of impacting
  mechanical oscillators.
\newblock {\em Nonlinearity}, 14:1517--1542, 2001.

\bibitem{GoKu05b}
W.~Govaerts, Yu.A. Kuznetsov, and A.~Dhooge.
\newblock Numerical continuation of bifurcations of limit cycles in {M}atlab.
\newblock {\em SIAM J. Sci. Comput.}, 27(1):231--252, 2005.

\bibitem{DoCh07}
E.~Doedel, A.~Champneys, T.~Fairgrieve, Y.~Kuznetsov, B.~Oldeman,
  R.~Paffenroth, B.~Sandstede, X.~Wang, and C.~Zhang.
\newblock {\em AUTO-07P: Continuation and Bifurcation Software for Ordinary
  Differential Equations.}, 2007.
\newblock {\em http://indy.cs.concordia.ca/auto}.

\bibitem{Er02}
B.~Ermentrout.
\newblock {\em Simulating, Analyzing, and Animating Dynamical Systems: {A}
  Guide to {XPPAUT} for Researchers and Students.}
\newblock SIAM, Philadelphia, 2002.

\bibitem{DaSc13}
H.~Dankowicz and F.~Schilder.
\newblock {\em Recipes for Continuation.}
\newblock SIAM, Philadelphia, PA, 2013.

\bibitem{DaZh05}
H.~Dankowicz and X.~Zhao.
\newblock Local analysis of co-dimension-one and co-dimension-two grazing
  bifurcations in impact microactuators.
\newblock {\em Phys. D}, 202:238--257, 2005.

\bibitem{JiCh17}
H.~Jiang, A.S.E. Chong, Y.~Ueda, and M.~Wiercigroch.
\newblock Grazing-induced bifurcations in impact oscillators with elastic and
  rigid constraints.
\newblock {\em Int. J. Mech. Sci.}, 127:204--214, 2017.

\bibitem{MaPi09}
J.F. Mason and P.T. Piiroinen.
\newblock The effect of codimension-two bifurcations on the global dynamics of
  a gear model.
\newblock {\em SIAM J. Appl. Dyn. Syst.}, 8(4):1694--1711, 2009.

\bibitem{ThZh06}
P.~Thota, X.~Zhao, and H.~Dankowicz.
\newblock Co-dimension-two grazing bifurcations in single-degree-of-freedom
  impact oscillators.
\newblock {\em J. Comput. Nonlinear Dynam.}, 1(4):328--335, 2006.

\bibitem{YiWe20}
S.~Yin, G.~Wen, J.~Ji, and H.~Xu.
\newblock Novel two-parameter dynamics of impact oscillators near degenerate
  grazing points.
\newblock {\em Int. J. Non-Linear Mech.}, 120:103403, 2020.

\bibitem{ZhDa06}
X.~Zhao and H.~Dankowicz.
\newblock Unfolding degenerate grazing dynamics in impact actuators.
\newblock {\em Nonlinearity}, 19(2):399--418, 2006.

\bibitem{DiKo02}
M.~di~Bernardo, P.~Kowalczyk, and A.~Nordmark.
\newblock Bifurcations of dynamical systems with sliding: Derivation of
  normal-form mappings.
\newblock {\em Phys. D}, 170:175--205, 2002.

\bibitem{PaIn10}
E.~Pavlovskaia, J.~Ing, M.~Wiercigroch, and S.~Banerjee.
\newblock Complex dynamics of bilinear oscillator close to grazing.
\newblock {\em Int. J. Bifurcation Chaos}, 20(11):3801--3817, 2010.

\bibitem{DiBu08}
M.~di~Bernardo, C.J. Budd, A.R. Champneys, and P.~Kowalczyk.
\newblock {\em Piecewise-smooth Dynamical Systems. Theory and Applications.}
\newblock Springer-Verlag, New York, 2008.

\bibitem{FrNo97}
M.H. Fredriksson and A.B. Nordmark.
\newblock Bifurcations caused by grazing incidence in many degrees of freedom
  impact oscillators.
\newblock {\em Proc. R. Soc. A}, 453:1261--1276, 1997.

\bibitem{FrNo00}
M.H. Fredriksson and A.B. Nordmark.
\newblock On normal form calculation in impact oscillators.
\newblock {\em Proc. R. Soc. A}, 456:315--329, 2000.

\bibitem{No91}
A.B. Nordmark.
\newblock Non-periodic motion caused by grazing incidence in impact
  oscillators.
\newblock {\em J. Sound Vib.}, 145(2):279--297, 1991.

\bibitem{Me07}
J.D. Meiss.
\newblock {\em Differential Dynamical Systems.}
\newblock SIAM, Philadelphia, 2007.

\bibitem{FoBi94}
S.~Foale and S.R. Bishop.
\newblock Bifurcations in impact oscillations.
\newblock {\em Nonlin. Dyn.}, 6:285--299, 1994.

\bibitem{LaBu94}
H.~Lamba and C.J. Budd.
\newblock Scaling of {L}yapunov exponents at nonsmooth bifurcations.
\newblock {\em Phys. Rev. E}, 50(1):84--90, 1994.

\bibitem{BaIn09}
S.~Banerjee, J.~Ing, E.~Pavlovskaia, M.~Wiercigroch, and R.K. Reddy.
\newblock Invisible grazings and dangerous bifurcations in impacting systems:
  The problem of narrow-band chaos.
\newblock {\em Phys. Rev. E}, 79:037201, 2009.

\bibitem{SiAv20}
D.J.W. Simpson, V.~Avrutin, and S.~Banerjee.
\newblock {N}ordmark map and the problem of large-amplitude chaos in impact
  oscillators.
\newblock {\em Phys. Rev. E}, 102:022211, 2020.

\bibitem{AlGe03}
E.L. Allgower and K.~Georg.
\newblock {\em Introduction to Numerical Continuation Methods.}
\newblock SIAM, Philadelphia, 2003.

\bibitem{Do07b}
E.~Doedel.
\newblock Lecture notes on numerical analysis of nonlinear equations.
\newblock In B.~Krauskopf, H.M. Osinga, and J.~Gal\'{a}n-Vioque, editors, {\em
  Numerical Continuation Methods for Dynamical Systems.}, pages 1--50.
  Springer, New York, 2007.

\bibitem{Go00b}
W.~Govaerts.
\newblock {\em Numerical Methods for Bifurcations of Dynamical Equilibria.}
\newblock SIAM, Philadelphia, 2000.

\end{thebibliography}

\end{document}